\documentclass[12pt]{article}
\usepackage{amsmath}
\usepackage{amsfonts}
\usepackage{amssymb}
\usepackage{latexsym}
\usepackage{amscd}
\usepackage{color}

\input amssym.def
\newcommand{\norma}[1]{\| #1 \|}
\newcommand{\conj}[2]{\left \{ {#1} \, : \, {#2} \right \}}

\newcommand{\cvf}{\overset{\omega}{\rightarrow}}

\newcommand {\cvfe} {\overset{\omega^\ast}{\rightarrow}}

\newcommand {\R}{\mathbb{R}}

\newcommand {\N} {\mathbb{N}}

\newcommand{\sol}[1]{\text{sol}(#1)}

\newcommand{\proof}{\noindent \textit{Proof: }}
\newcommand{\qed}{\hfill \ensuremath{\Box}}

\newtheorem{df}{Definition}[section]

\newtheorem{prop}[df]{Proposition}
\newtheorem{teo}[df]{Theorem}
\newtheorem{cor}[df]{Corollary}
\newtheorem{lemma}[df]{Lemma}
\newtheorem{ex}[df]{Example}

\newtheorem{rema}[df]{Remark}

\begin{document}

\title{Grothendieck-type subsets of Banach lattices}
\author{Pablo Galindo\footnote{Partially supported by Spanish MINECO/FEDER PGC2018-094431-B-I00.} ~and V. C. C. Miranda \footnote{Supported by a CAPES PhD's scholarship.}}
\date{}

\maketitle

\begin{abstract}
{In the setting of Banach lattices the weak (resp. positive) Grothendieck spaces have been defined. We localize such notions by defining  new classes of sets that we study and compare with some quite related different classes. This allows us to introduce and compare the corresponding  linear operators.

\noindent \textbf{Keywords: } almost Grothendieck sets, almost Grothendieck operators, disjoint operators, positive Grothendieck property,  positive Grothendieck sets, weak Grothendieck property

\noindent \textbf{ Mathematics Subject Classification (2010) } 46B42, 47B65.

\noindent
}\end{abstract}

\section{Introduction and preliminaries}
Throughout this paper $X$ and $Y$ will denote Banach spaces, $E$ and $F$ will denote Banach lattices. We denote by $B_X$ the closed unit ball of $X$.   Recall that $X$ is a Grothendieck space if weak* null sequences in $X'$ are weakly null.

 In a Banach lattice, two  further Grothendieck properties have  been considered. A Banach lattice $E$ has the  \textit{weak Grothendieck property} (resp. \textit{positive Grothendieck property}) if every disjoint weak* null sequence in $E'$ is weakly null (resp. every positive weak* null sequence in $E'$ is weakly null ).
 We refer to \cite{wnuk} and \cite{mach} for definitions and some results concerning those two properties.

A subset $A$ of a Banach space $X$ is said to be \textit{Grothendieck} if $T(A)$ is relatively weakly
compact in $c_0$ for every bounded operator $T: X \to c_0$. It is known that $X$ has the Grothendieck property if and only if $B_X$ is a Grothendieck set. Of course the Grothendieck property implies both the positive and the weak Grothendieck properties.
  The lattice $c$ of all real convergent sequences has the positive Grothendieck but it fails to have the weak Grothendieck property. On the other hand, $\ell_1$ is a Banach lattice with the weak Grothendieck property without the positive Grothendieck.

Keeping this $c_0$-valued operators point of view, we introduce and study a new class of sets in Banach lattices- that we name  \textit{almost Grothendieck} (see Definition \ref{ags1})-  and which  characterizes the weak Grothendieck property. In an analogous way, the notion of \textit{positive Grothendieck} set is defined.

Lattice versions of  Dunford-Pettis sets and  limited sets arise in Banach lattices theory in a natural way.
Following \cite{bouras} (resp. \cite{chen}),
a bounded subset $A \subset E$ is called \textit{almost Dunford-Pettis}
(resp. \textit{almost limited}) if every disjoint weakly null sequence in $E'$
converges uniformly to zero on $A$ (resp. if every disjoint weak* null sequence in $E'$ converges uniformly to zero on $A$). Obviously every almost limited set is almost Dunford-Pettis. It follows easily from the definitions that if $E$ has the weak Grothendieck property, then every almost Dunford-Pettis set in $E$ is almost limited.

We discuss also the relationship of the almost Grothendieck sets with these previously studied classes of subsets. Noticeably, all them can be characterized according to their behaviour against $c_0$-valued operators. This explains why we were led to some of the results included in the last sections.

We refer the reader  to \cite{alip, meyer} for background on Banach lattices theory and positive operators.
\section{Almost Grothendieck Sets}

Every bounded linear  operator $T: E \to c_0$ is uniquely determined by a weak* null sequence $(x_n') \subset E'$ such that $T(x) = \big(x_n'(x)\big)$ for all $x \in E.$ When this  sequence is disjoint in the dual Banach lattice $E',$ we  say that $T$ is a \textit{disjoint operator}.

\begin{df} \label{ags1}
We say that $A \subset E$ is an \textit{almost Grothendieck set} if $T(A)$ is relatively weakly compact in $c_0$ for every disjoint operator $T: E \to c_0$.
\end{df}

It is obvious that every Grothendieck set in a Banach lattice is almost Grothendieck. Obviously, each almost Grothendieck subset of $c_0$ is relatively weakly compact.

The class of almost Grothendieck sets satisfies the so-called Grothendieck \textit{encapsulating} property, i.e. if for every $\varepsilon>0$ there is an almost Grothendieck set $A_\varepsilon \subset E$ such that $A\subset A_\varepsilon +\varepsilon B_E,$ then $A$ is an almost Grothendieck set. This holds because the relatively weakly compact sets do have the Grothendieck encapsulating property \cite[Lemma 2, p. 227]{diestel}.

From the very definitions, we have the following characterization:

\begin{prop} \label{ags2}
For a Banach lattice $E$, the following are equivalent:
\begin{enumerate}
    \item $E$ has the weak Grothendieck property.

    \item Every disjoint operator $T: E \to c_0$  is weakly compact.

    \item $B_E$ is an almost Grothendieck set.
\end{enumerate}
\end{prop}

\proof
$(1) \Rightarrow (2)$ Let $T: E \to c_0$ be a disjoint operator with $T(x) = (x_n'(x))$ for all $x$. Since $(x_n')$ is a disjoint weak* null sequence in $E'$, we have that $(x_n')$ is weakly null. Then $T$ is weakly compact.

$(2) \Rightarrow (3)$ Obvious.

$(3) \Rightarrow (1)$ Let $(x_n')$ be a disjoint weak* null sequence in $E'$ and consider the disjoint operator $T(x) = (x_n'(x))$. Since $B_E$ is almost Grothendieck, $T(B_E)$ is a relatively weakly compact set in $c_0$. Therefore $x_n' \cvf 0$.
\qed

As a consequence of Proposition \ref{ags2}, $B_{\ell_1}$ is an almost Grothendieck set which is neither relatively weakly compact nor Grothendieck.

It is known that $A \subset X$ is a limited set if and only if $T(A)$ is relatively compact in $c_0$ for every bounded operator $T: X \to c_0$.
By using the concept of disjoint operators, we obtain a similar characterization for almost limited sets.
\begin{prop} \label{ags3}
A bounded set $A \subset E$ is almost limited if and only if $T(A)$ is relatively compact in $c_0$ for every disjoint operator $T: E \to c_0$.
\end{prop}

\proof
Assume that $A$ is almost limited in $E$ and let $T: E \to c_0$ be a disjoint operator, i.e. there exists a disjoint weak* null sequence $(x_n') \subset E'$ such that $T(x) = (x_n'(x))$ for all $x$.
For each $n$, let
$$ s_n = \sup\conj{|b_n|}{b = (b_1, b_2, \dots) \in T(A)} = \sup \conj{|x_n'(x)|}{x \in A} = \norma{x_n'}_A. $$
Since $A$ is almost limited, $\norma{x_n'}_A \to 0$, hence $T(A)$ is relatively compact in $c_0$, see \cite[Ex. 14, p. 168]{alip}.

Now, assume that $T(A)$ is relatively compact in $c_0$ for every disjoint operator from $E$ into $c_0$. Let $(x_n') \subset E'$ be a disjoint weak* null sequence and let $T: E \to c_0$ be the disjoint operator given by $T(x) = (x_n'(x))$ for all $x \in E$. Since $T(A)$ is relatively compact in $c_0$, by \cite[Ex. 14, p. 168]{alip} we have that
$\sup \conj{|b_n|}{b = (b_1, b_2, \dots) \in T(A)} \to 0$ as $n \to \infty$. Then $\norma{x_n'}_A \to 0$. Since $(x_n')$ is an arbitrary disjoint weak* null sequence in $E'$, $A$ is almost limited.
\qed

As a consequence of Proposition \ref{ags3}, every almost limited set is almost Grothendieck. On the other hand, $B_{\ell_2}$ is an almost Grothendieck set which is not almost limited. Note that the class of almost Grothendieck sets does not have any direct relation with the class of almost Dunford-Pettis sets. For instance, $B_{\ell_1}$ is an almost Grothendieck set which is not almost Dunford-Pettis, and $B_{c_0}$ is an almost Dunford-Pettis which is not almost Grothendieck.

\begin{rema} \label{rem1} If $T:E\to c_0$ is a positive operator and $(x_m)_m$ is a disjoint almost limited sequence, then $\{T(x_m)\}$ is a relatively compact subset of $c_0.$
\end{rema}

\proof Put $x'_n(x):=T(x)_n \,, x\in E.$ Suppose $\lim_n \sup_m\{|x_n'(x_m)|\}\neq 0.$ WLOG we may assume that there is $\varepsilon>0$ such that $\sup_{m}\{|x_n'(x_m)|\}>\varepsilon$ for all $n.$ So there is $m_n$ such that $|x_n'(x_{m_n})|>\varepsilon.$
We apply \cite[Ex. 22 p. 77]{alip} to obtain
a   disjoint sequence $(g_n)\subset E^*$ such that $0\leq g_n \leq x_n'$ and $g_n(x_{m_n})=x_n'(x_{m_n}).$ Further $(g_n)$ is a weak*-null sequence, thus it must converge uniformly to $0$ on the almost limited set $\{x_m \}.$ A fact prevented by the estimate $|g_n(x_{m_n})|>\varepsilon.$ Consequently, $\lim_n \sup_m\{|x_n'(x_m)|\}= 0,$ which means that $\{T(x_m): m\in \mathbb{N}\}$ is relatively compact.
\qed
\smallskip

If $E$ is $\sigma$-Dedekind complete, then we do not need to assume that $T$ is positive in Remark \ref{rem1}. Using the same notation, we  have that $x_n' \cvfe 0$ in $E'$ and $(g_n) \subset E'$ is disjoint with $|g_n| \leq |x_n'|$ for all $n$. It follows from Lemma 2.2 of \cite{chen} that $g_n \cvfe 0$. Then we proceed as in Remark \ref{rem1}. In other words, if $E$ is $\sigma$-Dedekind complete, then every disjoint almost limited sequence is a limited set.

\smallskip

Notice that the analogue of this remark for almost Grothendieck sequences fails even if we just would claim that the image is relatively weakly compact: Consider the sequence of unit vectors $(e_m)_m \subset B_{\ell_1}$ and the mapping $T:\ell_1 \to c_0$ given by $T\big((\alpha_i)_i\big)=\big((\sum _{i=j}^\infty \alpha_j)_j\big).$ Then $T(e_m)=(1,\stackrel{m}\dots,1,0,\dots)$ and clearly $\{T(e_m): m\in \mathbb{N}\}$ is not a  relatively weakly compact subset of $c_0.$

\smallskip

 Next we study how almost Grothendieck sets behave in the sum of two Banach lattices.
\begin{prop} \label{ags4}
Let $E$ and $F$ be two Banach lattices.
\begin{enumerate}
    \item Then $A \subset E$ and
$B \subset F$ are almost Grothendieck sets of $E$ and $F$, respectively if and only if $A \times B$ is an almost Grothendieck set in $E \oplus F$.

    \item Let $C \subset E \oplus F$ be an almost Grothendieck set. If $\pi_E: E \oplus F \to E$ and $\pi_F: E \oplus F \to F$ are the canonical projections, then $\pi_E(C)$ and $\pi_F(C)$ are almost Grothendieck in $E$ and $F$, respectively.
\end{enumerate}
\end{prop}

\proof
(1) We first assume that $A \subset E$ and $B \subset F$ are almost Grothendieck sets.
Let $T: E \oplus F \to c_0$ be a disjoint operator and let $(\varphi_n) \subset (E\oplus F)'$ be the disjoint weak* null sequence such that $T(x,y) = (\varphi_n(x,y))_n$ for all $(x,y) \in E \times F$. Since $(E \oplus F)' = E' \oplus F'$, there exist two weak*null disjoint sequences $(x_n') \subset E'$ and $(y_n') \subset F'$ such that $\varphi_n(x,y) = x_n'(x) + y_n'(y)$ for all $(x,y) \in E \oplus F$ and for each $n$.
Let $R: E \to c_0$ and $S: F \to c_0$ be the disjoints operators given by $R(x) = (x_n'(x))$ for every $x \in E$ and $S(y) = (y_n'(y))$ for every $y \in F$. Since $A$ and $B$ are almost Grothendieck sets, it follows that $R(A)$ and $S(B)$ are both relatively weakly compact sets in $c_0$. Hence $T(A) = R(A) + S(B)$ is a relatively weakly compact set in $c_0$.

Assume that $A \times B \subset E \oplus F$ is an almost Grothendieck set. Let $T: E \to c_0$ be a disjoint operator and let $(x_n') \subset E'$ be the disjoint weak* null sequence associated to $T$. In particular, $\varphi_n(x,y) = x_n'(x)$ defines a disjoint weak* null sequence in $(E \oplus F)'$. Hence we can consider the disjoint operator $R: E \oplus F \to c_0$ associated to $(\varphi_n)$. Since $R(x,y) = T(x)$, we have that $T(A) = R(A \times B)$ is relatively weakly compact in $c_0$. Therefore $A$ is almost Grothendieck in $E$. The same argument shows that $B$ is almost Grothendieck in $F$.

(2)  Let $C \subset E \oplus F$ be an almost Grothendieck set. We prove that $\pi_E(C)$ is almost Grothendieck in $E$. Let $T: E \to c_0$ be a disjoint operator and let $(x_n')$ be the disjoint weak* null sequence in $E'$ such that $T(x) = (x_n'(x))$ for all $x \in E$. It suffices us to prove that $T \circ \pi_E : E \oplus F \to c_0$ is a disjoint operator. We first note that
\begin{align*}
    |x_n' \circ \pi_E| (x,y) & = \sup \conj{|x_n' \circ \pi_E (z,w)|}{|(z,w)| \leq (x,y)} \\
        & = \sup \conj{|x_n'(z)|}{|z| \leq x} \\
        & = |x_n'|(x)
\end{align*}
holds for all $n$ and every $x \in E^+$ and $y \in F^+$. Now, if $(x,y) \in (E\oplus F)^+$ and $n\neq m$,
\begin{align*}
    [|x_n' \circ \pi_E| \wedge |x_m' \circ \pi_F|] (x,y) & = \inf \{ |x_n' \circ \pi_E| (z_1,w_1) + |x_m' \circ \pi_F| (z_2, w_2): \\
            & \quad \quad \quad  (z_i, w_i) \in (E \oplus F)^+, \\
            & \quad \quad \quad (z_1,w_1) + (z_2, w_2) = (x,y)\} \\
        & = \inf \conj{|x_n'|(z_1) + |x_m'|(z_2)}{z_i \in E^+, \, z_1 + z_2 = x} \\
        & = |x_n'| \wedge |x_m'| (x) = 0
\end{align*}
Therefore $(x_n' \circ \pi_E)$ is a disjoint sequence in $(E \oplus F)'$. Thus $\pi_E(C)$ is an almost Grothendieck set in $E$ as we claimed. By the same argument we can see that $\pi_F(C)$ is almost Grothendieck in $F$. \qed
\smallskip

In \cite{mach}, the authors proved that every L-space has the weak Grothendieck property. Also in \cite{miranda2}, it is proved that
the Banach lattice $(\bigoplus_{n \in \N} \ell_2^n)$ has the weak Grothendieck property. By using the same argument as in Proposition 2.8 of \cite{miranda2}, we have the following:

\begin{prop} \label{ags5}
Assume that $E'$, the dual Banach lattice, has an order unit, i.e. there exists $e' > 0$ such that $B_{E'} = [-e',e']$. Then $E$ has the weak Grothendieck property.
\end{prop}


Next proposition is a small improvement of Proposition 4.10 in \cite{mach}.

\begin{prop} \label{ags6}
If $E'$, the dual of some Banach lattice $E$, does not contain any copy of $\ell_1$, then $E$ has the weak Grothendieck property.
\end{prop}

\proof
Let $(x_n') \subset E'$ be a disjoint weak* null sequence and suppose it is not weakly null. Then, WLOG, we can assume that there exists $x'' \in E''$ such that $|x''(x_n')| \geq \epsilon$ for all $n$ form some $\epsilon > 0$. Since $E$ does not contain any copy of $\ell_1$, we can assume that $(x_n')$ is weakly Cauchy in $E'$. As $B_{E'''}$ is a $w^\ast$-compact subset of $E'''$, we have that there exists $x''' \in E'''$ such that $x'_n \cvfe x'''$ in $E'''.$ Moreover, since $(x_n')$ is a disjoint  sequence in $E'$, $(x_n')$ also is a disjoint sequence in $E'''$. Since $E'''$ is Dedekind complete, it follows that $x''' = 0$  by Proposition 1.4 in \cite{wnuk}. \qed

\smallskip

The converse of Proposition \ref{ags6} does not hold. For example, $\ell_\infty'$ contains a (complemented) copy of $\ell_1$ and $\ell_\infty$ has the weak Grothendieck property.

\begin{teo}
A Banach lattice $E$ has the weak Grothendieck property if and only if there does not exist a disjoint surjection $T: E \to c_0$.
\end{teo}

\proof
Assume that $E$ has the weak Grothendieck property. If $T: E \to c_0$ would be a disjoint surjection, $T(B_E)$ would be an open set by the Open Mapping Theorem and also a relatively weakly compact set by Proposition \ref{ags2}. All this  would imply $c_0$ is  a reflexive space.

%

Assume, by  way of contradiction, that there exists a disjoint weak* null sequence $(x_n') \subset E'$ which is not weakly null. If $(x_n')$ would be weakly Cauchy, we could argue as in the proof of above Proposition \ref{ags6} to conclude that $(x_n') \subset E'$ would be weakly null. Thus $(x_n')$ is not weakly Cauchy and by Roshental's $\ell_1$-theorem, we get that $(x_n')$ has a subsequence, that will also be denoted by $(x_n')$, that is equivalent to the canonical sequence $(e_n) \subset \ell_1$. Then, there is an isomorphism $S: \ell_1 \to E'$ such that $S e_n = x_n'$. But if $T: E \to c_0$ is the disjoint operator given by $T(x) = (x_n'(x))$ for all $x \in E$, then $T' = S$. By Exercise 2.49 in \cite{fabian}, we get that $T$ is a surjection. A contradiction.
\qed

\medskip

 Now we discuss whether the solid hull of an almost Grothendieck set is an almost Grothendieck set. This fails in general as the following simple example points out.

\begin{ex} \label{ags8}
Since $c$ does not have the weak Grothendieck property, $B_c = \sol{\{e\}}$ is not an almost Grothendieck set. However, the singleton $\{ e \}$ is almost Grothendieck.
\end{ex}


\begin{lemma} \label{ags9}
Let $T: E \to c_0$ be a disjoint order bounded operator such that $T(x) = (x_n'(x))$. Then  $|T|: E \to c_0$ is also disjoint and $|T|(x) = (|x_n'|(x))$ holds for all $n$.
\end{lemma}

\proof
Let $(x_n') \subset E'$ be the disjoint weak* null sequence such that $Tx = (x_n'(x))$ for all $x \in E$.
Since $c_0$ is Dedekind complete, there exists $|T|$ (see Theorem 1.13 in \cite{alip}). Note that if $B \subset c_0$ is such that there exists $(b_n) = \sup B$, then $b_m = \sup \conj{x_m}{x = (x_n) \in B}$ for all $m$. As $|T|(x) = \sup \conj{|Ty|}{|y| \leq x}$ holds for all $x \geq 0$, we have that
$$ \pi_m(|T|(x)) = \sup \conj{\pi_m(|Ty|)}{|y| \leq x} = \sup \conj{|x_n'(y)|}{|y| \leq x} = |x_n'|(x). $$
Therefore $|T|(x) = (|x_n'|(x))$ is a disjoint operator.
\qed

\smallskip

A Banach lattice $E$ has the \textit{property (d}) if $|x_n'| \cvfe 0$ for every disjoint weak* null sequence $(x_n')$. This property was introduced by Elbour in \cite{elbour}. It is easy to see that every Banach lattice with the weak Grothendieck property has the property (d).
Theorem 2.7 in \cite{miranda} shows that $E$ has the property (d) if and only if every disjoint linear operator on $E$ is regular, i.e. it can be written as the sum of two positive operators.

\begin{teo} \label{ags10}
Let $E$ be a Banach lattice with the property (d) and let $A \subset E$. Then  $|A| = \conj{|x|}{x \in A}$ is almost Grothendieck if and only if $\sol{A}$ is also almost Grothendieck.
\end{teo}

\proof
Since subsets of almost Grothendieck sets also are almost Grothendieck sets and $|A| \subset \sol{A}$, we have the sufficient condition.

Now assume that $|A|$ is almost Grothendieck and
let $T: E \to c_0$ be a disjoint operator. By Theorem 2.7 of \cite{miranda}, $T$ is a regular, hence, order bounded operator. Lemma \ref{ags9} yields that $|T|$ is also a disjoint operator. Besides, $T(\sol{A}) \subset \sol{|T|(|A|)}$ [Indeed, if $x \in \sol{A}$, there exists $y \in A$ such that $|x| \leq |y|$. So $|T(x)| \leq |T|(|x|) \leq |T|(|y|)$ what implies that $T(x) \in \sol{|T|(|A|)}$].
 Thus $|T|(|A|)$ is a relatively weakly compact subset of $(c_0)^+.$ Since $c_0$ has order continuous norm, which implies that $c_0$ is an ideal in $c_0''$, we may apply Theorem 4.39 in \cite{alip}, to conclude that $\sol{|T|(|A|)}$ is also relatively weakly compact. Consequently, $\sol{A}$ is almost Grothendieck.
\qed
\smallskip

It follows immediately from Theorem \ref{ags10} that if $A \subset E^+$ is an almost Grothendieck set in a Banach lattice with the property (d) then $\sol{A}$ is also an almost Grothendieck set.

\begin{rema} \label{ags11}
Let $E$ be a Banach lattice satisfying the following properties:
\begin{enumerate}
    \item $E$ has order continuous norm.
    \item Every almost Grothendieck set is relatively weakly compact.
    \item The lattice operations in $E$ are weakly sequentially continuous, i.e. $x_n \cvf 0$ in $E$ implies that $|x_n| \cvf 0$ in $E$.
\end{enumerate}
If $A \subset E$ is an almost Grothendieck set, then $\sol{A}$ is almost Grothendieck.
\end{rema}

\proof
If $A \subset E$ is an almost Grothendieck set, then $A$ is relatively weakly compact. Since the lattice operations in $E$ are weakly sequentially continuous, then $|A|$ is relatively weakly compact as well. On the other hand, as $E$ has order continuous norm, $E$ is an ideal in $E''$. By Theorem 4.39 in \cite{alip}, we have that $\sol{A} = \sol{|A|}$ is relatively weakly compact, hence almost Grothendieck.
\qed

\smallskip

In particular, $c_0$ satisfies the hypothesis in the above remark. In addition, it is obvious that if $F$ has the weak Grothendieck property, then the solid hull of every almost Grothendieck is almost Grothendieck.

 \smallskip


\begin{prop} \label{ags12}
If $I$ is an infinite index set, then $c_0(I)$  satisfies the hypothesis in Remark \ref{ags11}
\end{prop}
\proof
We first observe, despite it is well-known, that $c_0(I)$ has order continuous norm. Let $x'' \in (c_0(I))'' = \ell_\infty(I)$ and $y \in c_0(I)$ such that $|x''| \leq |y|$. If $x'' = (x_i'')$ and $y = (y_i)$, we have that $|x_i''| \leq |y_i|$ for all $i \in I$. Thus $(x_i'')$ vanishes at infinity. Therefore $c_0(I)$ is an ideal in $(c_0(I))''$, what implies that $c_0(I)$ has order continuous norm.

Let $A \subset c_0(I)$ be an almost Grothendieck set and let $(x_n) \subset A$. We want to prove that $(x_n)$ has a weakly convergent subsequence in $c_0(I)$. For each $n$, there exist a countable subset $J_n \subset I$ such that $x_{n,j} = 0$ for all $j \in I \setminus J_n$. If $J = \bigcup_{n \in \N} J_n$, then $J$ is a countable set. So we can find a bijection $F: \N \to J$. For each $j$, let $x_j'(x) = x_j$. Then, $(x_j') \subset (c_0(I))' = \ell_1(I)$ is a disjoint weak* null sequence. Define the disjoint operator $T(x) = (x_{F(k)}'(x))_k$. Since $(x_n)$ is an almost Grothendieck sequence, WLOG, we can assume that there exists $a \in c_0$ such that $Tx_n \cvf a$ in $c_0$.
Let $y = (y_j)_{j \in I}$ be defined by $y_j = 0$ if $j \notin J$ and $y_j = a_{F^{-1}(j)}$ if $j \in J$. Now if $x' = (x_j') \in c_0(I)$, then
\begin{align*}
    x'(x_n - y) & = \sum_{i \in I}                      x_i'(x_{n,i} - y_i) \\
                & = \sum_{j \in J} x_j'(x_{n,j} - y_j) \\
                & = \sum_{k \in \N} x_{F(k)}'(x_{n, F(k)} - a_k) \to 0.
\end{align*}

Finally, since $c_0(I)$ is an M-space, the lattice operations in $c_0(I)$ are sequentially weakly continuous \cite[Theorem 4.31]{alip}.
\qed
\begin{prop}\label{sol1} If $E$ satisfies the hypothesis in Remark \ref{ags11} and  $F$ is a Banach lattice with the weak Grothendieck property, then the solid hull of every almost Grothendieck set in $E \oplus F$ is  almost Grothendieck. \end{prop}
 \proof
 This follows from the fact that if $C \subset E \oplus F$ is almost Grothendieck, by Proposition \ref{ags4}, we have that both $A = \pi_E(C)$ and $B = \pi_F(C)$ are almost Grothendieck sets in $E$ and $F$, respectively. Then also are $\sol{A}$ and $\sol{B}$. Another application of Proposition \ref{ags4} implies that $\sol{A} \times \sol{B}$ is an almost Grothendieck set in $E \oplus F$. Since $\sol{C} \subset \sol{A} \times \sol{B}$, we are done. [Indeed, if $(x,y) \in \sol{C}$, there exists $(z,w) \in C$ such that $|(x,y)| \leq |(z,w)|$ in $E \oplus F$. This implies that $|x| \leq |z|$ in $E$ and $|y| \leq |w|$ in $F$. Note that $w = \pi_E(z,w) \in A$ and $z = \pi_F(z,w) \in B$. Then $x \in \sol{A}$ and $y \in \sol{B}$]. \qed
\medskip

Now, we seek to obtain a similar result to Theorem \ref{ags10} concerning Grothendieck sets in Banach lattices.

\begin{prop} \label{ags16}
Let $E$ be a Banach lattice such that the lattice operations in $E'$ are weak* sequentially continuous, i.e. $x_n' \cvfe 0$ in $E'$ implies that $|x_n'| \cvfe 0$ in $E'$. Then $|A|$ is a Gro\-then\-dieck set if and only if $\sol{A}$ is a Gro\-then\-dieck set.
\end{prop}

\proof
We prove the nontrivial implication. Assume that $|A|$ is a Grothendieck set and let $T: E \to c_0$ be a bounded linear operator. Since the lattice operations in $E'$ are weak* sequentially continuous, $T$ is a regular operator [The same proof used in \cite[Theorem 2.7]{miranda} holds in this case]. Now we can proceed with the same argument as in Theorem \ref{ags10} in order to get that $\sol{A}$ is a Grothendieck set.
\qed

The following remark gives conditions so that the solid hull of a Grothen\-dieck set in a Banach lattice is a Grothendieck set as well.

\begin{rema} \label{ags17}
Let $E$ be a Banach lattice satisfying the following properties:
\begin{enumerate}
    \item $E$ is a KB-space.

    \item Every Grothendieck set in $E$ is relatively weakly compact.
\end{enumerate}
If $A \subset E$ is a Grothendieck set, then $\sol{A}$ is Grothendieck.
\end{rema}

\proof
Since $E$ is a KB-space, $E$ is a band in $E''$ \cite[Theorem 4.60]{alip}. Now, by Theorem 4.39 in \cite{alip}, every relatively weakly compact subset of $E$ has a relatively weakly compact solid hull. If $A \subset E$ is a Grothendieck set, it must be by assumption relatively weakly compact, hence $\sol{A}$ is relatively weakly compact as well, what implies that $\sol{A}$ is a Grothendieck set.
\qed
\smallskip

 A Banach space is said to have the \textit{weak Gelfand-Phillips property} (ab. wGP) if every Grothendieck set is relatively weakly compact \cite{leung}. In this paper 
Leung proved that every separable Banach space, has the wGP property. Consequently, every separable KB-space satisfies the conditions in Remark \ref{ags17}, e.g. $\ell_1$ and $\ell_1(\ell_2^n)$.


\smallskip

By shifting from disjoint operators to positive operators we are led to the following definition clearly suggested by the positive Grothendieck lattice notion.

 \begin{df} \label{ags13} A subset $A \subset E$ is said to be a \textit{positive Grothendieck set} if $T(A)$ is relatively weakly compact for every positive operator $T: E \to c_0$. \end{df} The following statements are immediate consequences:
\begin{enumerate}
    \item The class of positive Grothendieck sets in $c_0$ coincides with the class of the relatively weakly compact sets. Hence it also coincides with the class of the almost Grothendieck sets in $c_0.$
    \item $E$ has the positive Grothendieck property iff $B_E$ is a positive Grothen\-dieck set.
    \item If $E$ has an order unit, then $E$ has the positive Grothendieck property. Indeed, let $e \in E$ be the order unit and let $(f_n) \subset (E')^+$ be a weak* null sequence. In particular,
    $$ \norma{f_n} = \sup_{x \in B_E} |f_n(x)| = \sup_{x \in [-e,e]} |f_n(x)| = f_n(e) \to 0. $$
    Consequently, $\norma{f_n} \to 0$, hence $f_n \cvf 0$. This is stronger than the positive Grothendieck property it was  proved  in \cite{wnuk} that if $E$ has an order unit, then it has the Dual Positive Schur property .
    \item If $E$ is a Banach lattice with and order unit and the property (d), then $E$ has the weak Grothendieck property. This follows from above item and \cite[Proposition 4.9]{mach}.
\end{enumerate}
\begin{ex}\label{+Gc_0(I)} The positive Grothendieck sets in $c_0(I)$ are relatively weakly compact.
This follows in the same way as in the course of the proof of Proposition \ref{ags12}.\end{ex}

Now we see that every positive operator (and consequently every regular operator) preserves positive Grothendieck sets.

\begin{prop} \label{ags14}
If $T: E \to F$ is a positive operator and $A \subset E$ is a positive Grothendieck set, then $T(A)$ is a positive Grothendieck set in $F$.
\end{prop}
\proof
If $S: F \to c_0$ is positive, then $S \circ T : E \to c_0$ is also a positive operator since $T \geq 0$. As $A$ is positive Grothendieck, $ST(A)$ is relatively weakly compact in $c_0$. Consequently, $T(A)$ is a positive Grothendieck set in $F$.
\qed

\smallskip

We conclude this section by observing that the same argument from Theorem \ref{ags10} can be used to show the following:

\begin{prop} \label{ags15}
Let $A$ be a subset of a Banach lattice $E$. Then $|A|$ is a positive Grothendieck set if and only if $\sol{A}$ is also positive Grothendieck.
\end{prop}

\section{Almost Grothendieck Operators}

Recall that a bounded operator $T: X \to Y$ is said to be \textit{Grothendieck} if $T(B_X)$ is a Grothendieck set in $Y$. Or, equivalently, $T'y_n' \cvf 0$ in $X'$ for every weak* null sequence $(y_n') \subset Y'$. In a natural way, we introduce the class of almost Grothendieck operators.

\begin{df}
 \label{agop1}
A bounded operator $T: X \to F$ is said to be \textit{almost Grothen\-dieck} if $T'y_n' \cvf 0$ in $X'$ for every disjoint weak* null sequence $(y_n') \subset F'$.
\end{df}

It is immediate that every Grothendieck operator $T: X \to F$ is almost Grothendieck. The converse does not hold though. For example, the identity map $I_{\ell_1} : \ell_1 \to \ell_1$ is almost Grothendieck but not Grothendieck.

The class of almost Grothendieck operators enjoys the following ideal property:

\begin{prop} \label{agop2}

\begin{enumerate}
 \item The set of almost Grothendieck operators is a closed subspace of $L(E;F)$ that is closed under composition on the right.

    \item $R \circ T: E \to G$ is almost Grothendieck for every bounded operator $R: F \to G$ whose adjoint $R'$ preserves disjoint sequences and every almost Grothendieck operator $T: E \to F.$

\end{enumerate}
\end{prop}

\proof
\begin{enumerate}

\item If $(T_k)\subset L(E,F)$ is a sequence of almost Grothendieck operators that converges to $T\in L(E,F),$ then given $\varepsilon>0,$ we find $k$ such that $T(B_E)\subset T_k(B_E) + \epsilon B_F.$ Since $T_k(B_E)$ is an almost Grothendieck set, $T(B_E)$ is as well an almost Grothendieck set because of the encapsulating property.

     Let $T \circ S: G \to F$ be an almost Grothendieck operator and $S: G \to E$  a  bounded operator.
    If $(y_n') \subset F'$ is a disjoint weak* null sequence, then $T'y_n' \cvf 0$ in $F'$. Thus $S'T'y_n' \cvf 0$ in $G'$, i.e. $(T \circ S)' y_n' \cvf 0$.

    \item If $(y_n') \subset G'$ is a disjoint weak* null sequence, then $(R'y_n')$ is a disjoint weak* null sequence in $F'$. Hence $(R \circ T)'y_n' = T'R'y_n' \cvf 0$.
    \qed
\end{enumerate}

Now, we state the connection between almost Grothendieck operators and almost Grothendieck sets.
\begin{teo} \label{agop3}
For a bounded operator $T: X \to F$, the following are equivalent:
\begin{enumerate}
    \item $T$ is an almost Grothendieck operator.

    \item $S \circ T: X \to c_0$ is weakly compact for every disjoint operator $S: F \to c_0$.

    \item $T(A)$ is almost Grothendieck in $F$ for every bounded set $A \subset X$.

    \item $T(B_X)$ is almost Grothendieck in $F$.
\end{enumerate}
\end{teo}

\proof
$(1) \Rightarrow (2)$ Let $S: F \to c_0$ be a disjoint operator and write $Sy = (y_n'(y))$ for all $y \in F$. Note that $S\circ T(x) = (y_n'(Tx)) = (T'y_n'(x))$. Since $T$ is almost Grothendieck, $T'y_n' \cvf 0$ in $X'$, thus $S \circ T$ is a weakly compact operator.

$(2) \Rightarrow (3)$ If $S: F \to c_0$ is a disjoint operator, by hypothesis, $S \circ T$ is weakly compact. Then $S(T(A))$ is relatively weakly compact for every bounded set $A \subset X$. Hence $T(A)$ is almost Grothendieck.

$(3) \Rightarrow (4)$ Obvious.

$(4) \Rightarrow (1)$ If $(y_n')$ is a disjoint weak* null sequence in $F'$, consider the disjoint operator $S: F \to c_0$ such that $Sy = (y_n'(y))$ for all $y \in F$. Since $T(B_X)$ is almost Grothendieck, $S(T(B_X))$ is relatively weakly compact. Then $S \circ T : X \to c_0$ is a weakly compact operator. On the other hand, since $S \circ T(x) = (T'y_n'(x))$, it follows that $T'y_n' \cvf 0$ in $X'$.
\qed

\smallskip

\begin{cor} \label{agop4}
The following assertions are equivalent:
\begin{enumerate}
    \item $F$ has the weak Grothendieck property.

    \item For all Banach space $X$, every bounded operator $T: X \to F$ is almost Grothendieck.

    \item For all Banach lattice $E$, every bounded operator $T: E \to F$ is almost Grothendieck.

    \item For all Banach lattice $E$, every positive operator $T: E \to F$ is almost Grothendieck.

    \item The identity operator $I_F$ is almost Grothendieck.
\end{enumerate}
\end{cor}

\proof
$(1) \Rightarrow (2)$ Let $T: X \to F$ be a bounded operator. If $S: F \to c_0$ is disjoint, since $F$ has the almost Grothendieck property, $S$ is a weakly compact operator. Then $S \circ T$ is also weakly compact. By Theorem \ref{agop3}, $T$ is almost Grothendieck.

$(2) \Rightarrow (3) \Rightarrow (4) \Rightarrow (5)$ Obvious.

$(5) \Rightarrow (1)$ If $I_F: F \to F$ is almost Grothendieck, then $B_F$ is almost Grothendieck. Proposition \ref{ags5} yields that $F$ has the weak Grothendieck property.
\qed

\smallskip

Theorem \ref{agop3} implies that every weakly compact operator is almost Gro\-then\-dieck. The converse does not hold. For instance, the identity map $i: \ell_1 \to \ell_1$ is not weakly compact.
It is natural to wonder whether the class of almost Grothendieck operators coincides with the class of weakly compact operators.
In order to answer this question, we need a lemma:
\begin{lemma} \label{agop6}
Let $(x_n) \subset E$ be an almost (resp. positive) Grothendieck sequence and let $T: \ell_1 \to E$ be the bounded operator given by $Ta = \sum_{j=1}^\infty a_j x_j$ for all $a = (a_j) \in \ell_1$. Then $T(B_{\ell_1}) $ is an almost (resp. positive) Grothendieck set.
\end{lemma}
\proof
If $S: F \to c_0$ is a disjoint operator, then $A = \conj{Sx_n}{n \in \N}$ is relatively weakly compact in $c_0$. Note that
\begin{align*}
    S(T(B_{\ell_1})) & = \conj{S(Ta)}{a = (a_j) \in B_{\ell_1}} \\
        & = \conj{S(\sum_{j=1}^\infty a_j x_j)}{\sum_{j=1}^\infty |a_j| \leq 1} \\
        & = \conj{\sum_{j=1}^\infty a_j Sx_j}{ \sum_{j=1}^\infty |a_j| \leq 1}
\end{align*}
which is the closed convex circled hull of $A$. By Krein-Smulian's theorem
it follows that $T(B_{\ell_1})$ is almost Grothendieck.
\qed

\begin{teo} \label{agop7}
For a Banach lattice $E$, the following assertions are equivalent:
\begin{enumerate}
    \item Every almost Grothendieck subset of $E$ is relatively weakly compact.

    \item For all Banach spaces $X$, every almost Grothendieck operator $T: X \to E$ is weakly compact.

    \item Every almost Grothendieck operator $T: \ell_1 \to E$ is weakly compact.
\end{enumerate}
\end{teo}
\proof
$(1) \Rightarrow (2)$
Assume that every almost Grothendieck subset of $E$ is relatively weakly compact. If $T: X \to E$ is an almost Grothendieck operator, then $T(B_E)$ is almost Grothendieck, hence relatively weakly compact. Thus $T$ is weakly compact.

$(2) \Rightarrow (3)$ Obvious.

$(3) \Rightarrow (1)$
If $A \subset E$ is an almost Grothendieck which is not relatively weakly compact, then we can find a sequence $(x_n) \subset A$ which has no weakly convergent subsequence. Let $T: \ell_1 \to E$ be defined by $T(a) = \sum_{j=1}^\infty a_j x_j$ for all $a = (a_j) \in \ell_1$.
Since $(x_n)$ is an almost Grothendieck sequence, Lemma \ref{agop6} yields that $T$ is almost Grothendieck. However, $T$ is not weakly compact because $(e_n) \subset \ell_1$ is bounded and $Te_n = x_n$ for all $n$.
\qed

\smallskip

Note that $c_0(I)$ and reflexive Banach lattices $E$ satisfy the properties in Theorem \ref{agop7} and they are also satisfied by  $c_0 \oplus E.$  It is clear that if a Banach lattice satisfies the hypothesis of Theorem \ref{agop7}, then it has the wGP property. As a separable space, $\ell_1$ has the wGP property, but it does not satisfy
the conditions in Theorem \ref{agop7}.

\smallskip

In the class of the  positive linear operators in  Banach lattices, we have a dominated type problem. For instance, let $S,T: E \to F$ be  positive operators such that $S \leq T$. The question is, if $T$ has some property $(\ast)$, does $S$ also have it?
We are interested in ascertaining what happens if $T$ is an almost Grothendieck operator.
\begin{prop} \label{agop9}
Let  $0 \leq S \leq T: E \to F$ with $T$ almost Grothendieck. If $F$ has the property (d), then $S$ also is almost Grothendieck.
\end{prop}
\proof
We observe that
$0 \leq S \leq T$ implies $0 \leq S' \leq T'$. If $y_n' \cvfe 0$ is disjoint in $F'$ then $|y_n'| \cvfe 0$. Thus $T' \,|y_n'| \cvf 0$ in $E'$. We claim that $S'y_n' \cvf 0$ in $E'$. Indeed, if $x' \in E'$, then
$$ |x'(S'y_n')|  \leq |x'| \, (S' \, |y_n'|) \leq |x'| \, (T' \,|y_n'|) \to 0.$$
\qed

 Recall that Grothendieck sets in Banach spaces are preserved by bounded linear operators. This fact does not happen for almost Grothendieck sets in Banach lattices. The  positive operator $T: \ell_1 \to c_0$ mentioned before Proposition \ref{ags4}  is not relatively weakly compact. Thus $T(B_{\ell_1})$ is not almost Grothendieck in $c_0$, in spite of  $B_{\ell_1}$ being an almost Grothendieck set.

  One may also  wonder under which conditions or which classes of bounded operators preserve almost Grothendieck sets. It is immediate that almost Grothendieck operators and disjoint operators preserve almost Grothendieck sets. Moreover, it is trivial that if $E$ is a Banach lattice that satisfies the  conditions of Theorem \ref{agop7}, then every bounded linear operator on $E$ into any other Banach lattice preserves almost Grothendieck sets.

\begin{prop} \label{agop10}
Let $T: E \to F$ be a linear operator such that $T': F' \to E'$ takes disjoint sequences onto disjoint sequences. If $A \subset E$ is almost Grothendieck, then also is $T(A)$.
\end{prop}

\proof
Let $S: F \to c_0$ be a disjoint operator. Let $(y_n') \subset F'$ be the disjoint weak* null sequence such that $S(y) = (y_n'(y))_n$ for all $y \in F$. Since $T'$ takes disjoint sequences onto disjoint sequences, $(T'y_n')_n$ is a disjoint weak* null sequence in $E'$. Then $R(x) = (T'y_n'(x))_n$ defines a disjoint operator from $E$ into $c_0$. Thus
$$ S(T(A)) = \conj{(y_n'(T(x)))}{x \in A} = R(A) $$
is relatively weakly compact.
\qed

\begin{cor} \label{agop5}
If $E$ has the weak Grothendieck property, then every bounded operator $T: E \to F$ whose adjoint preserves disjoint sequences is almost Grothendieck.
\end{cor}

\proof
By Corollary \ref{agop4}, the identity map $I_E$ is almost Grothendieck. On the other hand, since $T': F' \to E'$ preserves disjoint sequences, it follows from Proposition \ref{agop2} that $T = T \circ I_E: E \to F$ is almost Grothendieck.
\qed

\smallskip

Recall that $T: E \to F$ is called \textit{interval preserving} (resp. \textit{Riesz homomorphism}) whenever $T$ is positive and $T[0,x] = [0, T(x)]$ for all $x \in E^+$ (resp. $T(x \vee y = T(x) \vee T(y)$ holds for all $x,y \in E$). Also, if $T$ is an interval preserving operator, then $T'$ is a Riesz homomorphism. As Riesz homomorphisms preserve disjoint sequences, it follows that interval preserving operators satisfy the hypothesis in Proposition \ref{agop10}.
If $E\subset F$ is an (order) ideal, then the inclusion mapping is  interval preserving, and as a consequence of the above comments, an almost Grothendieck set in $E$ is also an almost Grothendieck set in $F.$ This is the case of $E$ a Banach lattice with order continuous norm and $F=E^{''}.$

\begin{prop} \label{agop11}
If the positive operator $S: E \to c_0$ is dominated by a disjoint operator $T: E \to c_0$, then $S$ is disjoint as well. Hence $S$ maps almost Grothendieck sets into relatively weakly compact ones.
\end{prop}

\proof
Write $T(x) = (x_n'(x))$ and $S(x) = (y_n'(x))$ for all $x \in E$. By assumption, $0 \leq y_n' \leq x_n'$ in $E'$ for all $n$. Thus $(y_n')$ is a disjoint sequence.
If $A$ is an almost Grothendieck set, then $S(A)$ is relatively weakly compact by definition.
\qed

\smallskip

Recall that a bounded operator $T: E \to F$ is said to be \textit{L-weakly compact} (resp. \textit{M-weakly compact}) if $T(B_E)$ is an L-weakly compact set in $F$, i.e. $\norma{y_n} \to 0$ for every disjoint sequence in $\sol{T(B_E)}$ (resp. if $\norma{Tx_n} \to 0$ for every disjoint sequence $(x_n) \subset B_E$).

\begin{prop} \label{agop12}
If $T: E \to F$ is an L-weakly compact operator, then $T'$ is almost Grothendieck.
\end{prop}
\proof
By Proposition 3.6.10 in \cite{meyer}, we have that $T': F' \to E'$ is M-weakly compact. If $(y_n')$ is a disjoint weak* null sequence in $E'$, then $(y_n')$ is a bounded disjoint sequence, what implies that $T'y_n' \to 0$ in $E'$. In particular, $T'y_n' \cvf 0$ in $E'$, and consequently, $T'$ is almost Grothendieck.
\qed
\medskip

We conclude this section with a number of results stressing the analogies between almost and positive Grothendieck sets. Next Proposition can be proved with the same argument used in \ref{agop7}

\begin{prop}\label{car+G}
Every positive Grothendieck subset of $E$ is relatively weakly compact if and only  every bounded linear operator $T: \ell_1 \to E$ such that $T(B_{\ell_1})$ is a positive Grothendieck set in $E$ is weakly compact.
\end{prop}

\begin{prop} \label{G+1}
Every positive Grothendieck set in $E$ is a Grothendieck set if and only if every bounded linear operator $T: X \to E$ such that $T(B_X)$ is a positive Grothendieck set in $E$ must be a Grothendieck operator, for all Banach space $X$.
\end{prop}
\proof 
We prove the nontrivial implication by contradiction. Let $A \subset E$ be a positive Grothendieck set which is not Grothendieck, then there exists a bounded linear operator $S: E \to c_0$ such that $S(A)$ is not relatively weakly compact. We note that, since $A$ is positive Grothendieck, $S$ cannot be a positive operator. In particular, we can find a sequence $(x_n) \subset A$ such that $(Sx_n)$ does not have any weakly convergent subsequence, and as a consequence $(x_n)$ cannot have any wakly convergent subsequence. By Lemma 6.2, we have that if $T: \ell_1 \to E$ is defined by $T(a) = \sum_{j=1}^\infty a_j x_j$, $a \in \ell_1$, then $T(B_{\ell_1})$ is a positive Grothendieck set in $E$. However, $T$ is not a Grothendieck operator; a contradiction.
\qed

\begin{prop}
Let $0 \leq S \leq T: E \to F$ be two bounded linear operators such that $T(B_{E})$ is a positive Grothendieck set in $F$. Then $S(B_E)$ is also a positive Grothendieck set in $F$.
\end{prop}

\proof
Since $T(B_E)$ is a positive Grothendieck set, $T'y_n' \cvf 0$ in $F'$ for every positive weak* null sequence $(y_n') \subset F'$. Hence $S'y_n' \cvf 0$ in $F'$ for every positive weak* null sequence  $(y_n') \subset F'$. Consequently, $S(B_E)$ is a positive Grothendieck set in $F$.
\qed

 \section{Polynomial considerations}

 Recent works studied polynomial versions of some classes of sets and properties in Banach lattices \cite{shiwangbu, wangshibu}. In this section, we are going to give a few contributions concerning polynomials and the properties studied in Section 2.
 Recall that a $k$-linear map $A: E_1 \times \cdots \times E_k \to F$ is said to be \textit{positive} if $A(x_1, \dots, x_k) \geq 0$ in $F$ for all $x_i \geq 0$ in $E_i$ for each $i = 1, \dots, k$. We say that a $k$-homogeneous polynomial $P: E \to F$ is positive if its associated symmetric $k$-linear map is positive. We refer  to \cite{loane, bubuskes12} for  background in polynomials on Riesz spaces/Banach lattices.

 We observe that positive polynomials do not preserve almost Grothendieck sets.
Let $p \in \N$ and let $P: \ell_p \to c_0$ be the $p$-homogeneous positive polynomial given by
 $ P ((x_i))=\left (\sum_{i=j}^\infty x_i^p\right)_j. $
 We note that $P(e_n)=(1, \stackrel{m}\dots, 1, 0, \dots)$ for all $n$ and  thus, the set $\{P(e_n)\}$ is not relatively weakly  compact in $c_0$. Consequently, $P(B_{\ell_p})$ is not almost Grothendieck in $c_0$.

Proposition \ref{ags14} shows that positive operators preserve positive Grothen\-dieck sets. A similar result does not hold for positive polynomials. For instance, if $P: \ell_2 \to c_0$ is the above $2$-homogeneous positive polynomial,  i.e., $P((x_i)) = \left ( \sum_{i=1}^\infty x_i^2 \right )_j$, then $P(B_{\ell_2})$ is not positive Grothendieck in $c_0$, but $B_{\ell_2}$ is positive Grothendieck in $\ell_2$.
\smallskip

We also note that a polynomial version of Remark \ref{rem1} holds.
\begin{prop} \label{rem2} If $P:E\to \mathbb{R}$ is a regular polynomial and  $(x_m)_m$ is a disjoint  almost limited weakly null sequence, then $\lim_m P(x_m)=0.$
\end{prop}
\proof
We first prove the result for positive polynomials.
It suffices to show the result for $k$-homogeneous polynomials. We argue by induction on $k.$ It is obvious for $k=1, $ so suppose it holds for $k-1$ and that $P\in P(^k E).$ Let $A$ denote the symmetric $k$-linear form associated with $P.$ Consider the linear forms on $E$ defined according to $f_m(x):=A(x,x_m, \dots, x_m).$ We check that $\lim_m f_m(x)=0$ by verifying  it for $x\geq 0$ and then by linearity. The polynomial $Q_x:E\to \mathbb{R}$ given by $Q_x(y)=A(x,y\dots, y)$ is a positive one, hence by the induction hypothesis, $0=\lim_m Q_x(x_m)=\lim_m A(x,x_m,\dots, x_m),$ as wanted. Consider now, the linear operator $T:E\to c_0$ given by $T(x):=\big((f_m(x))\big).$ It is a positive operator, thus by  Remark \ref{rem1}, the weakly null sequence $\{T(x_m)\}$ lies in a compact subset of $c_0.$ Hence $\lim_m \|T(x_m)\|=0,$ that is $\lim_m \sup_k |f_k(x_m)|=0,$ which leads to $\lim_m P(x_m)=0.$

If $P: E \to \R$ is a regular polynomial, then we can write $P = P_1 - P_2$, where $P_1, P_2: E \to \R$ are positive polynomials. Thus the result follows immediately.
\qed
\smallskip

In Proposition \ref{rem2}, we can not change almost limited sequences by almost Grothendieck sequences. Let $p \in \N$ with $p > 1$ and let $P: \ell_p \to c_0$ be the positive polynomial $P((x_i)) = \left ( \sum_{i=j}^\infty x_i^p \right )_j$. In particular, $(e_n) \subset \ell_p$ is a disjoint almost Grothendieck weakly null sequence and $(P(e_n))$ is not relatively weakly compact.

Gonz\'{a}lez and Guti\'{e}rrez studied polynomial versions of the Grothendieck property in a Banach space \cite{gongut}. Our next result,  whose proof  is inspired in the proof of Theorem 14 in \cite{gongut}, is a polynomial characterization of the positive Grothendieck property in a Banach lattice. For a polynomial $P$ on $E$ we denote by $\widetilde{P}$ its Davie-Gamelin extension to $E''$ \cite{DG}.

 \begin{teo} The following statements are equivalent.
 \begin{enumerate}
     \item $E$ has the positive Grothendieck property.

     \item for every integer $k$, given a positive sequence $(P_n) \subset \mathcal{P}^r(^k E)$ with $P_n(x) \to 0$ for all $x \in E$, then $\widetilde{P_n}(x'') \to 0$ for all $x'' \in E''$.

     \item for some integer $k$, given a positive sequence $(P_n) \subset \mathcal{P}^r(^k E)$ with $P_n(x) \to 0$ for all $x \in E$, then $\widetilde{P_n}(x'') \to 0$ for all $x'' \in E''$.
 \end{enumerate}
 \end{teo}

 \proof
 $ (1) \Rightarrow (2) $ For $k=1$, the result is obvious.
 Assume it holds for $k-1$ and  let us prove it for $k$. Let $(P_n)$ be a positive sequence in $\mathcal{P}^r(^k E)$ such that $P_n(x) \to 0$ for all $x \in E$. If $A_n$ is the symmetric $k$-linear form associated to $P_n$, then the polarization formula yields that
 $ A_n(x_1, \dots, x_n) \to 0 $ for all $x_1, \dots, x_n \in E$. Fixing $x_1 \in E^+$, we define $Q_n(x) = A_n(x_1, x, \dots, x)$. In particular, $(Q_n)$ is a sequence of positive  polynomials in $\mathcal{P}^r(^{k-1}E)$ such that $Q_n(x) \to 0$ for all $x \in E$. By  the induction hypothesis, we have that $\widetilde{Q_n}(x'') \to 0$ for all $x'' \in E''$.  Note that $B_n(x_2, \dots, x_k) = A_n(x_1, x_2, \dots, x_k)$ is the $(k-1)$-linear symmetric form associated to $Q_n$. Thus $$\widetilde{Q_n}(x'') = \widetilde{B_n}(x'', \dots, x'') = \widetilde{A_n} (x_1, x'', \dots, x''),$$  where $\widetilde{A_n}$ denotes the multilinear symmetric form associated to $\widetilde{P_n}$. Again by  the polarization formula,
 $ \widetilde{A_n}(x_1, x_2'', \dots, x_n'') \to 0 $ for all $x_2'', \dots, x_n'' \in  E''$ and fixed $x_1 \in E^+$.

 On the other hand, fixing $x_2'', \dots, x_k'' \in E''$, we can consider the sequence $(x_n')$ in $E'$ given by $x_n'(x) = \widetilde{A_n}(x, x_2'', \dots, x_n'')$. Since $(x_n')$ is positive and weak* null, we have   by assumption, that $x'_n \cvf 0$ in $E'$. Hence
 $$ \widetilde{A_n}(x_1'', \dots, x_n'') \to 0 $$
 for all $x_1'', \dots, x_n'' \in E''$.

 $(2) \Rightarrow (3)$ Obvious.

 $(3) \Rightarrow (1)$ Let $(x_n')$ be a positive weak* null sequence in $E'$. Then $P_n(x) = (x_n'(x))^k$, $x \in E$, defines a sequence of positive polynomials. Since $x_n'\cvfe 0$ in $E'$, $P_n(x) \cvf 0$ for all $x \in E$. Therefore $\widetilde{P_n}(x'') \to 0$ for all $x'' \in E''$. In particular, if $x'' \in E''$, then $(x''(x_n'))^k = \widetilde{P_n}(x'') \to 0$.
 \qed

\section{Miscelanea}

This section contains additional links of the notions we have discussed and related results.

Note that $B_{\ell_1}$ is an almost Grothendieck set which is not Grothendieck. Consequently, $I_{\ell_1}$ is an almost Grothendieck operator which is not a Grothen\-dieck operator.

The same idea used in the proof of Proposition \ref{G+1} can be used to obtain the following result.

\begin{prop}
For a Banach lattice $E$, the following statements are equivalent:
\begin{enumerate}
    \item Every almost Grothendieck set in $E$ is a Grothendieck set.

    \item For every Banach space $X$, every almost Grothendieck operator $T: X \to E$ is a Grothendieck operator.
\end{enumerate}
\end{prop}

\begin{ex}
Let $E$ be a Banach lattice with an order unit. If $E$ has the property (d), then every almost Grothendieck set subset of $E$ is almost limited. {\rm Indeed, since $E$ has an order unit, every norm bounded subset of $E$ is order bounded. Hence, every almost Grothendieck set is order bounded, so it lies in an order interval. Therefore, every almost Grothendieck set  is almost limited according to \cite[see Prop. 2.2]{miranda}.}
\end{ex}
Note that $\ell_\infty(\ell_2^n)$ and $E'$, where $E$ is an L-space, are Banach lattices with an order unit and property (d).
\smallskip

Recall that a bounded linear operator $T: X \to E$ is said to be \textit{almost limited} if $\norma{T' y_n'} \to 0$ for every disjoint weak* null sequence $(y_n') \subset E'$, or equivalently $T(B_X)$ is an almost limited set in $E$. It is obvious that every almost limited operator is almost Grothendieck. However, $I_{\ell_1}$ is not almost limited but it is almost Grothendieck. 

\begin{prop}
For a given Banach lattice $E$, the following statements are equivalent:
\begin{enumerate}
    \item Every almost Grothendieck set in $E$ is almost limited.

    \item For every Banach space $X$, every almost Grothendieck operator $T: X \to E$ is almost limited.

    \item Every almost Grothendieck operator $T: \ell_1 \to E$ is almost limited.
\end{enumerate}
\end{prop}

\proof
$(1) \Rightarrow (2) \Rightarrow (3)$ Obvious.

$(3) \Rightarrow (1)$ Let $A \subset E$ be an almost Grothendieck set. If $A$ is not almost limited, there exists a disjoint weak* null sequence $(x_n') \subset E'$ which does not converge uniformly to zero on $A$. WLOG, we can assume that there exist $\epsilon > 0$ and $(x_n) \subset A$ such that $|x_n'(x_n)| \geq \epsilon$ for all $n$. Let $T: \ell_1 \to E$ be given by $T(a) = \sum_{j=1}^\infty a_j \, x_j $ for all $a = (a_j) \in \ell_1$. By Lemma \ref{agop6}, $T$ is an almost Grothendieck operator. However, $T$ is not almost limited. Indeed, note that
$T': E' \to \ell_\infty = \ell_1'$ is such that
$$T'(x') (a) = x'(Ta) = x'(\sum_{j=1}^\infty a_j x_j) = \sum_{j=1}^\infty a_j x'(x_j).$$
Then
$T'(x') = (x'(x_j))_j$ for all $x' \in E'$. In particular, $(x_n') \subset E'$ is disjoint and weak* null, but $\norma{T'x_n'}_\infty \geq |x_n'(x_n)| \geq \epsilon$ for all $n$. \qed

\smallskip

 Recall that a Banach lattice $E$ is said to have the \textit{wDP property}  if every relatively weakly compact subset of $E$ is an almost DP set.
Also, recall that $E$ has the wDP property iff every weakly compact operator on $E$ is almost Dunford-Pettis, i.e. it maps disjoint weakly null sequences onto norm null ones \cite{sanchez}.

\begin{prop}
Let $E$ be a Banach lattice with the weak Grothendieck property. Then, $E$ has the wDP property if and only if every disjoint operator on $E$ is almost Dunford-Pettis.
\end{prop}

\proof
By Proposition \ref{ags2}, every disjoint operator on $E$ is weakly compact. Assuming that $E$ has the wDP property, if $T: E \to c_0$ is disjoint, then $T$ is almost Dunford-Pettis.
Now assume that every disjoint operator on $E$ is almost Dunford-Pettis. Let $(x_n) \subset E$ and $(x_n') \subset E'$ be disjoint weakly null sequences. Consider the disjoint operator on $E$ given by $T(x) = (x_n'(x))$. Since $x_n' \cvf 0$ in $E'$, $T$ is a weakly compact operator, what implies that $T$ is almost DP. Therefore, $\norma{Tx_n}_\infty \to 0$. Consequently, $x_n'(x_n) \to 0$. \qed

\smallskip

%
%
%
The following is similar to Remark \ref{rem1}.
\begin{prop} \label{adp1}
If $(x_m)_m$ is a disjoint almost Dunford-Pettis sequence in $E$, then it is a Dunford-Pettis set.
\end{prop}
\proof Let $T: E \to c_0$ be a weakly compact operator. We are going to check that $\{T(x_m)_m\}$ is a  relatively compact subset of $c_0.$ If $T(x) = (x_n'(x))$ for all $x \in E,$ we have that  $x_n' \cvf 0$ in $E'$. Suppose $\lim_n \sup_m \{ |x_n'(x_m)| \} \neq 0$. WLOG, we may assume that $\sup_m \{ |x_n'(x_m)| \} > \epsilon$ for some $\epsilon > 0$ and every $n$. So there exist $m_n$ such that $|x_n'(x_{m_n})| > \epsilon$.

Apply \cite[Ex. 22, p. 77]{alip} to obtain a disjoint sequence $(y_n') \subset E'$ such that $|y_n'| \leq |x_n'|$ and $y_n'(x_{m_n}) = x_m'(x_{m_n})$ for all $n$. Since $\{ x_n' \} \subset E'$ is relatively weakly compact and $(y_n') \subset \sol{ \{ x_n' \}}$ is disjoint, $y_n' \cvf 0$ in $E'$, see \cite[Theorem 4.34]{alip}. Therefore, $(y_n')$ must converge uniformly to zero on the almost Dunford-Pettis set $\{ x_m \}$, what contradicts the fact that $|y_n'(x_{m_n})| > \epsilon$.
\qed

\begin{cor}
 If $(x_m)_m$ is a disjoint almost Dunford-Pettis weakly null sequence and if $T: E \to Y$ is a weakly compact operator, then $\norma{Tx_m} \to 0$.
\end{cor}

\proof
Assume first that $Y$ is reflexive.
By using Proposition \ref{adp1}, $(x_m)$ is a Dunford-Pettis sequence in $E$. Consequently, $(Tx_m)\subset Y$ is a Dunford-Pettis  set, hence  limited in $Y$ because it is reflexive.
 In addition, reflexive Banach spaces have the GP property (i.e.  limited sets are relatively compact), therefore the weakly null sequence $(Tx_m)$ lies in a relatively compact set, hence it converges to $0.$

Since $T: E \to Y$ is weakly compact, there exists a reflexive Banach space $Z$ and two bounded linear operators $S: E \to Z$ and $R: Z \to Y$ such that $T = RS$ (Davis-Figiel-Johnson-Pelczynski's Theorem). By using the first part, $\norma{Sx_m} \to 0$. Consequently, $\norma{Tx_m} \leq \norma{R} \, \norma{Sx_m} \to 0$.
\qed

\smallskip

As we pointed after Remark \ref{rem1}, every disjoint almost limited sequence in a $\sigma$-Dedekind complete Banach lattice is a limited sequence as well. In Proposition \ref{adp1}, we saw that every disjoint almost Dunford-Pettis sequence in an arbitrary Banach lattice is a Dunford-Pettis sequence. We observe that, we do not have the analogue for disjoint almost Grothendieck sequences. For instance, the unit sequence $(e_m)_m \subset \ell_1$ is not Grothendieck. A similar fact neither holds for disjoint positive Grothendieck sequences, e.g. $(e_m) \subset c$ is not Grothendieck.

\smallskip

\noindent \textbf{Acknowledgments: } This paper is part of the second author Ph. D. thesis supervised by Prof. Mary Lilian Louren\c co whom he warmly thanks for her support and encouragement.





\noindent\emph{Contact:}

\noindent Pablo Galindo. Departamento de An\'{a}lisis Matem\'{a}tico. Universidad de Valencia. Spain. e-mail: galindo@uv.es

\smallskip

\noindent Vin\'icius C. C. Miranda, University of S\~ao Paulo, Brazil, e-mail: vinicius.colferai.miranda@usp.br

\end{document}